\documentclass[11pt]{amsart2000}
\usepackage{amssymb}
\newtheorem*{theorem}{Theorem}
\theoremstyle{definition}
\newtheorem*{observation}{Observation}
\numberwithin{equation}{section}

\begin{document}
\setlength{\unitlength}{0.01in}
\linethickness{0.01in}
\begin{center}
\begin{picture}(474,66)(0,0)
\multiput(0,66)(1,0){40}{\line(0,-1){24}}
\multiput(43,65)(1,-1){24}{\line(0,-1){40}}
\multiput(1,39)(1,-1){40}{\line(1,0){24}}
\multiput(70,2)(1,1){24}{\line(0,1){40}}
\multiput(72,0)(1,1){24}{\line(1,0){40}}
\multiput(97,66)(1,0){40}{\line(0,-1){40}}
\put(143,66){\makebox(0,0)[tl]{\footnotesize Proceedings of the Ninth Prague Topological Symposium}}
\put(143,50){\makebox(0,0)[tl]{\footnotesize Contributed papers from the symposium held in}}
\put(143,34){\makebox(0,0)[tl]{\footnotesize Prague, Czech Republic, August 19--25, 2001}}
\end{picture}
\end{center}
\vspace{0.25in}
\setcounter{page}{321}
\title{Quasiorders on topological categories}
\author{V\v{e}ra Trnkov\'a}
\address{Math. Institute of Charles University\\ 
Sokolovsk\'a 83\\
18675 Praha 8}
\email{trnkova@karlin.mff.cuni.cz}
\thanks{Financial support of the Grant Agency of the Czech Republic under
the grants no.201/99/0310 and 201/00/1466 is gratefully 
acknowledged. Supported also by MSM 113200007.}
\thanks{V\v{e}ra Trnkov\'a,
{\em Quasiorders on topological categories},
Proceedings of the Ninth Prague Topological Symposium, (Prague, 2001),
pp.~321--330, Topology Atlas, Toronto, 2002}
\begin{abstract} 
We prove that, for every cardinal number $\alpha\geq {\mathfrak c}$, there
exists a metrizable space $X$ with $|X|=\alpha$ such that for every pair
of quasiorders $\leq_1$, $\leq_2$ on a set $Q$ with $|Q| \leq \alpha$
satisfying the implication
$$q \leq_1 q' \implies q \leq_2 q'$$
there exists a~system $\{ X(q) : q\in Q\}$ of non-homeomorphic clopen
subsets of $X$ with the following properties: 
\begin{itemize}
\item $q \leq_1 q'$ if and only if $X(q)$ is homeomorphic to a~clopen
subset of $X(q')$,
\item $q \leq_2 q'$ implies that $X(q)$ is homeomorphic to a~closed
subset of $X(q')$ and
\item $\neg (q \leq_2 q')$ implies that there is no one-to-one
continuous map of $X(q)$ into $X(q')$.
\end{itemize}
\end{abstract}
\maketitle

\section{Introduction and the Main Result}

Let $\mathcal{M}$ be a~class of morphisms of a~category $\mathcal{K}$
containing all isomorphisms and closed with respect to the composition. 
Then $\mathcal{M}$ determines a~quasiorder $\preccurlyeq$ on the class of
objects of $\mathcal{K}$ by the rule
$$
\begin{array}{l}
X\preccurlyeq Y\ 
\mbox{if and only if there exists}\ 
m\in \mathcal{M}\\
\mbox{with the domain $X$ and the codomain $Y$.}
\end{array}
$$

An $\mathcal{M}${\it -representation of a~quasiordered set} $(Q,\leq )$ 
{\it in} $\mathcal{K}$ is any collection $X=\{ X(q) : q\in Q \}$ of
non-isomorphic objects of $\mathcal{K}$ such that, for every $q,q'\in Q$,
$$
q\leq q'
\text{ if and only if }
X(q)\preccurlyeq X(q').
$$

Which quasiordered sets have $\mathcal{M}$-representations in which
categories $\mathcal{K}$ for which classes $\mathcal{M}$? 
In topology, this investigation with $\mathcal{M}$ being the class of all 
homeomorphic embeddings has rather long tradition. 
In 1926, C. Kuratowski and W. Sierpi\'{n}ski proved (see \cite{3,4}) that
the antichain of the cardinality $2^{\mathfrak c}$ and the ordinal 
${\mathfrak c}^{+}$ have such representations within the category of all
subspaces of the real line ${\mathbb R}$. 
In 1993, P.T. Matthews and T.B.M. McMaster refreshed this field of 
problems and proved (see \cite{5}) that every partially ordered set of the
cardinality ${\mathfrak c}$ has such representation. 
In 1999, A.E. McCluskey, T.B.M. McMaster and W.S. Watson proved (see 
\cite{8}) that the set $(\exp {\mathbb R},\subseteq)$ of all subsets of 
the real line ${\mathbb R}$ ordered by the inclusion also has such a
representation, i.e.\ a~representation by subspaces of ${\mathbb R}$ with
respect to the embeddability. 
Since $(\exp {\mathbb R},\subseteq)$ contains an antichain of the
cardinality $2^{\mathfrak c}$, their result implies the first result of
C. Kuratowski and W. Sierpi\'{n}ski. 
Analogously it implies the previous result of P.T. Matthews and 
T.B.M. McMaster because $(\exp {\mathbb R},\subseteq)$ contains an
isomorphic copy of any partially ordered set of cardinality at most
${\mathfrak c}$. 
The authors of \cite{8} also announced that they have a~counterexample 
consistent with ZFC to the statement that every partially ordered set of 
the cardinality $2^{\mathfrak c}$ can be represented by subspaces of 
${\mathbb R}$ with respect to the embeddability. 
In \cite{7}, A.E. McCluskey and T.B.M. McMaster generalized the 
construction of \cite{8} and they proved that, for any infinite cardinal 
numbers ${\mathfrak a},{\mathfrak b}$ such that 
${\mathfrak b}^{{\mathfrak a}} = {\mathfrak b}$, every partially ordered
set $(P,\leq)$ with $|P| \leq {\mathfrak b}$ has a~representation (with
respect to the embeddability) by subsets of {\it every} $T_3$-space $X$
containing a~dense subset $D$ with $|D| \leq {\mathfrak a}$ whenever every
non-empty open set of $X$ has the cardinality ${\mathfrak b}$. 
Moreover, if $X$ admits a~homeomorphism into itself such that every 
$x\in X$ has an infinite orbit, then every quasiordered set obtained from 
a partially ordered set $(P,\leq)$ with $|P| \leq {\mathfrak b}$ by 
splitting every $p\in P$ into $2^{\mathfrak b}$ mutually comparable
points has also such a~representation.

Given a~cardinal number $\alpha$, A.E. McCluskey and T.B.M. McMaster also
present (in \cite{6}) a~construction of a~$T_0$-space $X$ with 
$|X| = \delta(\alpha)$ such that every quasiordered set $(Q,\leq)$ with
$|Q| \leq \alpha$ has a~representation by subsets of $X$ (with respect to
the embeddability), where $\delta (\alpha)$ is the smallest cardinal
number $\delta$ such that there exist $\alpha$ distinct cardinal numbers
$\gamma$ (not necessarily infinite) smaller than $\delta$. 
If $\alpha = \aleph_0$, then $\delta(\alpha) = \aleph_0$; hence there
exists a~countable $T_0$-space $X$ such that every countable quasiordered
set has a~representation (with respect to the embeddability) by subsets of
$X$.
For $\alpha$ uncountable, the size of $X$ is rather high. 
The statement below offers a~stronger result for any $\alpha \geq
{\mathfrak c}$.

For every cardinal number $\alpha \geq {\mathfrak c}$ there exists a
metrizable space $X$ with $|X| =\alpha$ such that every quasiordered set
$(Q,\leq )$ with $|Q| \leq \alpha$ has an $\mathcal{M}$-representation
by retracts of $X$ whenever
\begin{enumerate}
\item 
$\mathcal{M}$ consists of one-to-one continuous maps and it contains all
coproduct injections (= homeomorphic embeddings onto clopen subspaces) or
\item 
$\mathcal{M}$ consists of continuous surjections and it contains all 
product projections.
\end{enumerate}

The claim (1) implies e.g.\ the representability of every quasiordered set
$(Q,\leq)$ with $|Q| \leq \alpha$ by retracts of $X$ with respect to the
embeddability or the embeddability onto closed subspaces or onto retracts
or onto clopen subspaces. 
The claim (2) implies the representability e.g.\ by being continuous image
or a~quotient or a~continuous open image or a~factor in a~product. 
Here we prove a result much stronger than (1) above, namely the following.

\begin{theorem} 
For every cardinal number $\alpha\geq {\mathfrak c}$, there exists
a~metrizable space $X$ with $|X|=\alpha$ such that, for every pair
$\leq_1$, $\leq_2$ of quasiorders on a set $Q$ with $|Q| \leq \alpha$
satisfying the implication
$$q \leq_1 q' \implies q \leq_2 q',$$
there exists a~system $\{ X(q) : q\in Q\}$ of non-homeomorphic clopen
subsets of $X$ with the following properties:
\begin{enumerate}
\item 
$q \leq_1 q'$ if and only if $X(q)$ is homeomorphic to a~clopen subset of
$X(q')$,
\item 
$q \leq_2 q'$ implies that $X(q)$ is homeomorphic to a~closed subset of
$X(q')$ and
\item 
$\neg (q \leq_2 q')$ implies that there is no one-to-one continuous map of
$X(q)$ into $X(q')$.
\end{enumerate}
\end{theorem}

The proof of this theorem is presented in the part II of the present 
paper. 
Inspecting the proof, one can see that none of the spaces $X(q)$
has an isolated point. 
Hence if one adds $\alpha$ isolated points to every $X(q)$ (and also to
the space $X$), she (he) gets the following result:

{\em
For every cardinal number $\alpha \geq {\mathfrak c}$ there exists a
metrizable space $X$ with $|X| = \alpha$ such that every quasiordered set
$(Q,\leq)$ with $|Q| \leq \alpha$ has an $\mathcal{M}$-representation by
clopen subspaces of $X$ whenever $\mathcal{M}$ is the class of all
continuous bijections.
}

We do not present here the proof of the above claim (2) and of its 
stronger variant concerning simultaneous representation of a~pair
$\leq_1$, $\leq_2$ of quasiorders by product projections and continuous
surjections. 
The proof will appear elsewhere.

Finally, let us refresh here some results which have applications in the
present field of problems. 
In \cite{1}, J. Ad\'{a}mek and V. Koubek introduced a~sum-productive 
representation of a~partially ordered commutative semigroups as
follows: 
If $(S,\circ,\leq)$ is a~partially ordered commutative semigroup
(i.e.\ if $(S,\circ)$ is a~commutative semigroup and $\leq$ is a~partial
order on $S$ such that 
$$
(a\leq b)\
\mbox{and}\ 
(a'\leq b')
\implies 
a\circ a' \leq b\circ b'),$$
then its sum-productive representation in a~category $\mathcal{K}$ is any 
collection $\{ X(s) : s\in S\}$ of objects of $\mathcal{K}$
such that
\begin{itemize}
\item[($\times$)] 
$X(s\circ s')$ is always isomorphic to the product $X(s)\times X(s')$ and
\item[($\leq$)] 
$s\leq s'$ if and only if $X(s)$ is isomorphic to a~summand of $X(s')$
(i.e.\ $X(s')$ is a~coproduct of $X(s)$ and an object of $\mathcal{K}$).
\end{itemize}

In \cite{11}, J. Vin\'{a}rek proved that every partially ordered 
commutative semigroup has a~sum-productive representation in the category
of all metric zero-dimensional spaces. 
Every partially ordered set $(P,\leq)$ can be enlarged to a~partially 
ordered set $(S,\leq )$ in which every pair of elements has an infimum. 
Putting $p\circ p' = \inf\{ p,p'\}$, one gets the partially ordered 
commutative semigroup $(S,\circ,\leq)$ and the Vin\'{a}rek's result can be
applied. 
Hence

{\em
every partially ordered set has a representation by zero-dimensional
metrizable spaces with respect to the embeddability onto clopen subspaces.
}

In \cite{10}, sum-productive representations in the category of all
$F_{\sigma\delta}$-and-$G_{\sigma\delta}$ subspaces of the Cantor
discontinuum are examined. 
Omitting the product forming again, we get the result that

{\em 
if ${\mathbb C}$ is a~countable set, then $(\exp {\mathbb C},\subseteq)$
has a representation by $F_{\sigma\delta}$ and $G_{\sigma\delta}$ 
subspaces of the Cantor discontinuum with respect to the embeddability
onto clopen subspaces.
}

These scattered results are surrounded by many unsolved questions, such 
as: 
which {\it quasiordered\/} sets can be represented by metrizable 
zero-dimen\-sion\-al spaces or by closed or Borelian or all subsets of the 
real line with respect to the embeddability onto clopen subspaces or onto
Borelian subspaces or onto closed subspaces or onto retracts; 
and many others.

\section{Proof of the Theorem}

Let $\alpha\geq {\mathfrak c}$ be given. 
Let $Q$ be a~set with $| Q| =\alpha$. 
In the part A of the proof, we suppose that we have sets 
$S^{(1)},\ldots,S^{(4)}$ of metrizable spaces of cardinality $\alpha$
with the five properties below and we prove Theorem using such sets. 
In the part B, we prove that such sets $S^{(1)}, \ldots ,S^{(4)}$ really
do exist. 

\subsection*{Part A}

\subsubsection*{a}

Thus, let us suppose that $S^{(1)},\ldots,S^{(4)}$ are sets of metrizable
spaces with $|S^{(i)}|=\alpha$ and $|Y|=\alpha$ for every 
$Y\in S^{(i)}$, $i=1,\ldots,4$, such that the statements (1)--(5) below
are satisfied:
\begin{enumerate}
\item
$S^{(1)}=\{ A_q,B_q : q\in Q\}$ are spaces such that, for every 
$q\in Q$, $A_q$ is homeomorphic to a~clopen subspace of $B_q$ and $B_q$ is
homeomorphic to a~closed subspace but to no clopen subspace of $A_q$; 
moreover, if $q\neq q'$, there exists no continuous one-to-one map of any
$A_q,B_q$ into any $A_{q'},B_{q'}$;
\item
$S^{(2)}=\{ D_q,E_q : q\in Q\}$ are such that
$D_q$ is homeomorphic to a~closed subspace of $E_q$ and $E_q$ is
homeomorphic to a~closed subspace of $D_q$ but there exists no
homeomorphism of $D_q$ onto a~clopen subspace of $E_q$ and no
homeomorphism of $E_q$ onto a~clopen subspace of $D_q$; moreover, if
$q\neq q'$, then there exists no continuous one-to-one map of any
$D_q,E_q$ into any $D_{q'},E_{q'}$;
\item
$S^{(3)}=\{ M_q,N_q : q\in Q\}$ are spaces such that $M_q$ is
homeomorphic to a~clopen subspace of $N_q$ and $N_q$ is homeomorphic to
a~closed subspace but to no clopen subspace of $M_q$ (i.e.\ they are
mutually situated as $A_q$ and $B_q$ in (1)); 
moreover, if $q\neq q'$, there exists no continuous one-to-one map of any
$M_q,N_q$ into any $M_{q'},N_{q'}$;
\item
$S^{(4)}=\{ G_q,H_q : q\in Q\}$ are spaces such that $G_q$ is
homeomorphic
to a~clopen subspace of $H_q$ and $H_q$ is homeomorphic to a~clopen
subspace of $G_q$ but $G_q$ is not homeomorphic to $H_q$; 
moreover, if $q\neq q'$, there exists no continuous one-to-one map of any
$G_q,H_q$ into any $G_{q'},H_{q'}$;
\item
if $i,j=1,\ldots,4$, $i\neq j$, $Z\in S^{(i)}$, $Y\in S^{(j)}$,
then there exists no cont\-inuous one-to-one map of $Z$ into $Y$.
\end{enumerate}

We put 
$$X = 
\coprod\limits_{q\in Q}
\textstyle{
(A_q
\coprod B_q
\coprod (D_q\times \omega )
\coprod (E_q\times\omega )
\coprod M_q
\coprod N_q
\coprod G_q)
}
$$
where $\coprod$ and $\displaystyle\coprod$ denote the coproduct (= 
disjoint union as clopen subspaces) and $\omega$ is a countable discrete 
space (i.e.\ $D_q\times\omega$ is a coproduct of countably many copies of 
$D_q$ and analogously for $E_q$).

Hence $X$ is a~metrizable space and $|X|=\alpha$. 
We show that $X$ has all the required properties. 
In the reasoning below, we shall frequently use the statement (5) without
mentioning it. 

\subsubsection*{b}

Let $(Q,\leq_1,\leq_2)$ be a~set with two quasiorders $\leq_1,\leq_2$ such
that
$$q \leq_1 q' \implies q \leq_2 q'.$$
To begin with, let us suppose, moreover, that $\leq_1$ is a~partial order,
i.e.\
\begin{equation}\label{star}
(q \leq_1 q')\ \mbox{and}\ (q'\leq_1q) \implies q=q'
\end{equation}
(this requirement will be removed at the end of part A of the proof).

By means of $S^{(1)}$--$S^{(3)}$, we construct a system 
$\{ X(q) : q\in Q\}$ of clopen subspaces of $X$ such that
\begin{itemize}
\item $q\leq_1q'$ if and only if $X(q)$ is homeomorphic to a~clopen subset
of $X(q')$,
\item $q\leq_2q'$ if and only if $X(q)$ is homeomorphic to a~closed subset
of $X(q')$ and
\item if $\neg (q\leq_2q')$, then there is no continuous one-to-one map of
$X(q)$ into $X(q')$.
\end{itemize}
First, let us define $<$ by the rule
$$
q<q'\
\mbox{if and only if}\
(q\leq_1q')\
\mbox{and}\ 
(\neg (q'\leq_1q))\
\mbox{and}\
(q'\leq_2q),
$$
and $q\leq q'$ denotes $q<q'$ or $q=q'$. 
(By (\ref{star}), $\leq$ is a~partial order.) 
Let $\mathcal{C}$ be the system of all components of $\leq$ (i.e.\ 
$q,q'\in C\in \mathcal{C}$ if and only if there exist $q_0,\ldots,q_n$
in $Q$ such that $q=q_0\leq q_1 \geq q_2 \leq \ldots q_n = q'$).

For every $q\in C\in \mathcal{C}$ put
$$
\begin{array}{l}
S_q^{(1)} = 
\{ B_{q'} : q'\leq q\} \cup
\{ A_{q'} : q'\in C\text{ and }\neg (q'\leq q)\},\\
X_q^{(1)} = \coprod S_q^{(1)}.
\end{array}
$$
Moreover, let us denote
$$
\begin{array}{l}
B_C = \{ B_q : q\in C\},\\
D_C = \{ D_q\times\omega : q\in C\},\\
E_C = \{ E_q\times\omega : q\in C\}.
\end{array}
$$

\begin{observation} 
If $q\in C$ and $q'\in C'$ with $C,C'\in \mathcal{C}$, $C\neq C'$, then
there exists no continuous one-to-one map of $X_q^{(1)}$ into
$X_{q'}^{(1)}$. 
If $q,q'\in C$, then $q<q'$
if and only if $X_q^{(1)}$ is homeomorphic to a~clopen subspace of 
$X_{q'}^{(1)}$ and $X_{q'}^{(1)}$ is homeomorphic to a~closed but
not to any clopen subspace of $X_q^{(1)}$.
\end{observation}

\subsubsection*{c}

On $Q$, let us define $q \,R\, q'$ if and only if either
$$\begin{array}{l}
q < q'\ \mbox{or}\\
q \leq_2 q'\ \mbox{and}\ q' \leq_2 q\ 
\mbox{but neither}\ q \leq_1 q'\ \mbox{nor}\ q' \leq_1 q.
\end{array}$$
Let $\prec$ be the~transitive envelope of $R$, and let $\preceq$ mean
$\prec$ or $=$. 
Then $\preceq$ is a~quasiorder on $Q$; 
let $\mathcal{L}$ be the system of all its components.
Clearly, each $C\in \mathcal{C}$ is a~subset of precisely one 
$L\in \mathcal{L}$.
For every $q\in C\in \mathcal{C}$, $C\subseteq L$, put
$$
\begin{array}{l}
S_q^{(2)} =
S_q^{(1)} \cup
\bigcup\limits_{C'\in \mathcal{C}, C\neq C'\subseteq L}
(B_{C'}\cup E_{C'}\cup D_{C'}) \cup E_C;\\
X_q^{(2)} = \coprod S_q^{(2)}
\end{array}
$$

\begin{observation} 
If $q,q'\in C\in \mathcal{C}$, $q\neq q'$, then $S_q^{(2)}$ differs from
$S_q^{(1)}$ by the same summand 
$$
\bigcup\limits_{C'\in \mathcal{C}, C\neq C'\subseteq L}
(B_{C'}\cup E_{C'}\cup D_{C'})\cup E_C
$$ 
as $S_{q'}^{(2)}$ differs from $S_{q'}^{(1)}$. 
Hence $X_q^{(2)}$ and $X_{q'}^{(2)}$ are in the same relations as
$X_q^{(1)}$ and $X_{q'}^{(1)}$. 
If $C,C'\in \mathcal{C}$, $C\neq C'$ and both $C,C'$ are subsets of 
$L\in \mathcal{L}$, $q\in C$, $q'\in C'$, then $X_q$ is homeomorphic to a
closed subspace but not to any clopen subspace of $X_{q'}$ and vice
versa. 
This is because $\coprod (E_C\cup D_C\cup E_{C'})$ is homeomorphic to
a closed (but not to any clopen!) subspace of 
$\coprod (E_{C'}\cup D_{C'}\cup E_C)$, by (2) 
(and 
$\coprod (S_q^{(1)}\cup B_{C'})$ is homeomorphic to a~closed subspace of
$\coprod (S_{q'}^{(1)}\cup B_C)$, by (1)). 
\end{observation}

Hence $\{ X_q^{(2)} : q\in L\}$ is such that
$$\mbox{
$q \leq_1 q'$ if and only if $X_q^{(2)}$ is homeomorphic to a clopen
subset of $X_{q'}^{(2)}$
}$$
and
$$\mbox{
$q \leq_2 q'$ if and only if $X_q^{(2)}$ is homeomorphic to a closed
subset of $X_{q'}^{(2)}$.
}$$
(For every $q,q'\in L$, we have $q\leq_2q'$ and $q'\leq_2q$, of course.) 
Moreover, if $L,L'\in \mathcal{L}$, $L\neq L'$, $q\in L$, $q'\in L'$, then
there exists no one-to-one continuous map of $X(q)$ into $X(q')$ and vice
versa. 
Let us denote 
$$
S_L^{(2)} = \{ B_C\cup E_C\cup D_C : C\in \mathcal{C},C\subseteq L\}.
$$

\subsubsection*{d}

On $Q$, let us define $q \dot{R} q'$ if and only if either 
$$\begin{array}{l}
q\prec q'\ \mbox{or}\\ 
q \neq q', q \leq_1 q'\ \mbox{or}\\
q \neq q', q \leq_2 q'.
\end{array}$$
Let $\dot{\prec}$ be the~transitive envelope of $\dot{R}$, and let 
$\dot{\preceq}$ mean $\dot{\prec}$ or $=$. 
Let $\mathcal{T}$ be the system of all components of the quasiorder
$\dot{\preceq}$. 
Every $L\in\mathcal{L}$ is a~subset of precisely one $T\in \mathcal{T}$. 
If $L,L'\in \mathcal{L}$, $L\neq L'$, $q\in L$, $q'\in L'$ and
$q\leq_2q'$, then for no $\bar {q}\in L$, $\bar{q}'\in L'$ is
$\bar{q}' \leq_2 \bar{q}$ because $L\neq L'$. 
Hence in fact, $\dot{\preceq}$ determines a~partial order on the 
components in $\mathcal{L}$, let us denote it also by $\dot{\preceq}$. 
For every $L\in \mathcal{L}$ define 
$$
\widetilde S_q^{(2)} = 
\bigcup\limits_{L'\dot{\prec }L} \widetilde S_{L'}^{(2)}
$$
and for $q\in L$,
$$
\begin{array}{l}
\widetilde S_q^{(2)} = S_q^{(2)} \cup \widetilde S_L^{(2)},\\
\widetilde X_q^{(2)} = \coprod \widetilde S_q^{(2)}.
\end{array}
$$

If $q\in L,q'\in L'$ and $L,L'$ are components of $\mathcal{L}$
incomparable in the partial order $\dot{\preceq}$, then there is no
continuous one-to-one map of $\widetilde X_q^{(2)}$ into 
$\widetilde X_{q'}^{(2)}$ and vice versa. 
If $q,q'\in L$, then $\widetilde X_q^{(2)}$ is in the same position with
respect to $\widetilde X_{q'}^{(2)}$ as $X_q^{(2)}$ to $X_{q'}^{(2)}$. 
If $q'\in L' \dot{\preceq}L\ni q$, then $\widetilde X_{q'}^{(2)}$ is
homeomorphic to a~clopen subspace of $\widetilde X_q^{(2)}$ and there
exists no continuous one-to-one map of $\widetilde X_q^{(2)}$ into
$\widetilde X_{q'}^{(2)}$. 
This is what we need whenever $q'\leq_1q$. 
However, if $q'\leq_2q$ and $\neg (q'\leq_1q)$, we still have to modify
the systems $\widetilde S_q^{(2)}$ and the spaces $\widetilde X_q^{(2)}$
by means of the spaces in $S^{(3)}$ in (3). 
For $q\in L$, we put
$$
\begin{array}{l}
S_q^{(3)} = 
\widetilde S_q^{(2)} \cup 
\{ N_{q'} : q' \in T \text{ and } q' \leq_1 q\} \cup
\{ M_{q'} : q'\in T \text{ and } \neg (q'\leq_1q)\};\\
X(q) = \coprod S_q^{(3)}.
\end{array}
$$

First, we outline why these ``new summands'' do not destroy the mutual 
position of $\widetilde X_{q_1}^{(2)}$ and $\widetilde X_{q_2}^{(2)}$ 
whenever $q_1,q_2\in L\subseteq T$:
\begin{itemize}
\item
$q_1\leq_1q_2$: 
if $N_{q'}\in S_{q_1}^{(3)}$, i.e.\ $q' \leq_1q_1$, we get $q'\leq_1q_2$
so that $N_{q'}\in S_{q_2}^{(3)}$; 
since $S_{q_2}^{(3)}$ always contains either $M_{q'}$ or $N_{q'}$, the 
case $M_{q'}\in S_{q_1}^{(3)}$ is clear.
\item 
$q_1\leq_2q_2$: 
if $N_{q'}\in S_{q_1}^{(3)}$, i.e.\ $q' \leq_1q_1$, then either 
$q'\leq_1q_2$ hence $N_{q'}\in S_{q_2}^{(3)}$, or $\neg (q'\leq_1q_2)$
hence $M_{q'}\in S_{q_2}^{(3)}$; 
but $N_{q'}$ is homeomorphic to a~closed subspace of $M_{q'}$.
\end{itemize}
No other case for $q_1,q_2\in L$ is possible (if $\neg (q_1\leq q_2)$ for 
$q_1,q_2\in L$, then already $\widetilde S_{q_1}^{(2)}$ is not 
homeomorphic to a~clopen subspace of $\widetilde S_{q_2}^{(2)}$).

Next we show that the system $\{ X(q) : q\in Q\}$ has all the required
properties.

Let $q_1,q_2\in Q$. 
If $q_1,q_2\in L$, we have just proved it. 
If $q_1\in T,q_2\in T'$ where $T,T'$ are distinct elements of $\mathcal{T}$
(i.e.\ distinct components of the quasiorder $\dot{\preceq}$), then there
is no continuous one-to-one map of $X(q_1)$ into $X(q_2)$ or vice versa. 
But this precisely corresponds to the fact that $\neg (q_1\leq_2q_2)$ and
$\neg (q_2\leq_2q_1)$. 
The remaining case is that $q_1,q_2$ are in distinct components $L_1,L_2$
of $\mathcal{L}$ but both $L_1$ and $L_2$ are subsets of a~component 
$T\in \mathcal{T}$.
We discuss the following cases:
\begin{itemize}
\item[$\alpha$)] 
$L_1$ and $L_2$ are incomparable in the partial order $\dot{\preceq}$: 
then there exists no one-to-one continuous map of 
$\widetilde X_{q_1}^{(2)}$ into $\widetilde X_{q_2}^{(2)}$ and adding any
summands $M_{q'},N_{q'}$ cannot change this situation;
hence there exists no continuous one-to-one map of $X(q_1)$ into $X(q_2)$
or vice versa.
\item[$\beta$)] 
let $L_1\dot{\prec }L_2$: 
hence $\widetilde X_{q_1}^{(2)}$ is homeomorphic to a~clopen subspace of
$\widetilde X_{q_2}^{(2)}$ and there exists no one-to-one continuous map
of $\widetilde X_{q_2}^{(2)}$ into $\widetilde X_{q_1}^{(2)}$, adding of
any summands $M_{q'},N_{q'}$ cannot change the latter fact; since 
$\neg (q_2\leq_2q_1)$, this is precisely what we need. 
Thus we have to discuss the cases $q_1\leq_1q_2$ and $q_1\leq_2q_2$. 
However the reasoning is precisely the same as in the above discussed
cases when $q_1,q_2\in L$ and $q_1\leq_1q_2$ or $q_1\leq_2q_2$.
\end{itemize}

\subsubsection*{e}

Finally, we have to remove the condition (\ref{star}): 
If $\leq_1$ is a~quasiorder, we use the standard trick that we define an
equivalence on $Q$ by the rule
$$
q\sim q'
\text{ if and only if }
q\leq_1q'\text{ and }q'\leq_1q.
$$
Then $\leq_1$ determines a~partial order and $\leq_2$ a~quasiorder on the
set $Q/\!\!\sim$; 
let us denote them by $\leq_1$ and $\leq_2$ as well. 
For every $[q]\in Q/\!\!\sim$, we construct $S_{[q]}^{(3)}$ and 
$X_{[q]}^{(3)}=\coprod S_{[q]}^{(3)}$ as described and, for $q\in [q]$, we
put 
$$
S_q^{(3)} = 
S_{[q]}^{(3)} \cup \{ G_{q'} : q'\in Q\setminus \{ q\}\} \cup\{ H_q\}
$$
where $G_{q'},H_q$ are from $S^{(4)}$, i.e.\ they satisfy the condition (4). 
Then, for $X_q=\coprod S_q^{(3)}$, the system $\{ X_q : q\in Q\}$ has
all the required properties.

\subsection*{Part B}

\subsubsection*{}
To finish the proof, it remains to show that the systems 
$S^{(1)},\ldots,S^{(4)}$ with the properties (1)--(5) do exist. 
For this, we use the fact (see \cite{9}) that for every 
$\alpha \geq {\mathfrak c}$ there exists a~system ${\mathbb S}$ of
metrizable spaces of the cardinality $\alpha$ such that for every 
$Y,Z\in {\mathbb S}$ and every continuous map $f:Y\to Z$, either $f$ is
constant or $Y=Z$ and $f$ is the identity, and the system ${\mathbb S}$
itself is large enough: for our purpose it suffices
$|{\mathbb S}| = \alpha$ (although, as shown in \cite{9}, such 
${\mathbb S}$ with $|{\mathbb S}| = 2^{\alpha}$ does exist). 
Let
$$\{ A_q^{*},D_q^{*},E_q^{*},M_q^{*},G_q^{*} : q\in Q\}$$
be a~subsytem of ${\mathbb S}$.
\begin{enumerate}
\item
Choose two distinct points, say $a_q,b_q$, in the space $A_q^{*}$ and, in
the space $A_q^{*}\times\omega$, identify $(b_q,n)$ with $(a_q,n+1)$. 
The obtained space is $A_q$, then $B_q=A_q^{*} \coprod A_q$. 
Clearly, $S^{(1)}=\{ A_q,B_q : q\in Q\}$ satisfies (1).
\item
Choose two distinct points, say $d_q^{(1)},d_q^{(2)}$ in $D_q^{*}$ and
$e_q^{(1)},e_q^{(2)}$ in $E_q^{*}$; 
in $D_q^{*}\times \omega$, identify $(d_q^{(2)},n)$ with $(d_q^{(1)},n+1)$
and denote $\widetilde D_q$ the obtained space; analogously define
$\widetilde E_q$; then put
$$
\begin{array}{l}
D_q = D_q^{*}\coprod\widetilde D_q\coprod\widetilde E_q\\
E_q = E_q^{*}\coprod\widetilde D_q\coprod\widetilde E_q
\end{array}
$$
Then $S^{(2)}=\{ D_q,E_q : q\in Q\}$ satisfies (2).
\item
The spaces $M_q$ and $N_q$ are constructed from $M_q^{*}$ as $A_q$ and
$B_q$ from $A_q^{*}$ in (1). 
Then $S^{(3)}=\{ M_q,N_q : q\in Q\}$ satisfies (3).
\item
Let $K$ be a~closed subset of the Cantor discontinuum such that $K$ is 
homeomorphic to $K\coprod K\coprod K$ but not to $K\coprod K$. 
Such a~space was constructed in \cite{2}. 
Put $G_q=K\times G_q^{*}$ and $H_q=(K\coprod K)\times G_q^{*}$. 
Since $K$ is zero-dimensional while $G_q^{*}$ is connected and every 
continuous map of $G_q^{*}$ into itself is either the identity or a 
constant, every continuous one-to-one map $f:G_q\to H_q$ sends every
$K\times\{ a\}$ into $(K\coprod K)\times\{ a\}$ for an arbitrary 
$a\in G_q^{*}$. 
Since $K$ and $K\coprod K$ are non-homeomorphic, the space $G_q$ is not
homeomorphic to $H_q$. 
However $G_q$ is homeomorphic to a~clopen subspace of $H_q$ and vice
versa.
\item
The system $S^{(1)}\cup S^{(2)}\cup S^{(3)}\cup S^{(4)}$ satisfies (5) 
because
$$
\{ A_q^{*},D_q^{*},E_q^{*},M_q^{*},N_q^{*},G_q^{*} : q\in Q\}
\subseteq {\mathbb S},
$$
evidently.
\end{enumerate}


\begin{thebibliography}{10}

\bibitem{1}
J.~Ad{\'a}mek and V.~Koubek, \emph{On representations of ordered commutative
  semigroups}, Algebraic theory of semigroups (Proc. Sixth Algebraic Conf.,
  Szeged, 1976), North-Holland, Amsterdam, 1979, pp.~15--31. \MR{80h:06012}

\bibitem{3}
Kuratowski C., \emph{Sur la puissance de l'ensemble des ``nombres de
  dimension'' au sense de m. fr\`{e}chet}, Fund. Math. \textbf{8} (1926),
  201--208.

\bibitem{4}
Kuratowski C. and Sierpi\'{n}ski W., \emph{Sur un probl\`{e}me de m.fr\`{e}chet
  concernant les dimensions des ensembles lineaires}, Fund. Math. \textbf{8}
  (1926), 193--200.

\bibitem{2}
Jussi Ketonen, \emph{The structure of countable {B}oolean algebras}, Ann. of
  Math. (2) \textbf{108} (1978), no.~1, 41--89. \MR{58 \#10647}

\bibitem{5}
P.~T. Matthews and T.~B.~M. McMaster, \emph{Families of spaces having
  prescribed embeddability order-type}, Rend. Istit. Mat. Univ. Trieste
  \textbf{25} (1993), no.~1-2, 345--352 (1994). \MR{96e:54015}

\bibitem{7}
A.~E. McCluskey and T.~B.~M. McMaster, \emph{Realizing quasiordered sets by
  subspaces of ``continuum-like'' spaces}, Order \textbf{15} (1998/99), no.~2,
  143--149. \MR{2001i:54042}

\bibitem{6}
\bysame, \emph{Representing quasi-orders by embeddability ordering of families
  of topological spaces}, Proc. Amer. Math. Soc. \textbf{127} (1999), no.~5,
  1275--1279. \MR{99h:06002}

\bibitem{8}
A.~E. McCluskey, T.~B.~M. McMaster, and W.~S. Watson, \emph{Representing
  set-inclusion by embeddability (among the subspaces of the real line)},
  Topology Appl. \textbf{96} (1999), no.~1, 89--92. \MR{2000g:54070}

\bibitem{9}
Ale{\v{s}} Pultr and V{\v{e}}ra Trnkov{\'a}, \emph{Combinatorial, algebraic and
  topological representations of groups, semigroups and categories},
  North-Holland Publishing Co., Amsterdam, 1980. \MR{81d:18001}

\bibitem{10}
V{\v{e}}ra Trnkov{\'a}, \emph{Homeomorphisms of products of subsets of the
  {C}antor discontinuum}, Dissertationes Math. (Rozprawy Mat.) \textbf{268}
  (1988), 37. \MR{89e:54078}

\bibitem{11}
Ji{\v{r}}{\'\i} Vin{\'a}rek, \emph{Representations of commutative semigroups by
  products of metric $0$-dimensional spaces}, Comment. Math. Univ. Carolin.
  \textbf{23} (1982), no.~4, 715--726. \MR{84m:20072}

\end{thebibliography}
\providecommand{\bysame}{\leavevmode\hbox to3em{\hrulefill}\thinspace}
\providecommand{\MR}{\relax\ifhmode\unskip\space\fi MR }
\providecommand{\MRhref}[2]{%
  \href{http://www.ams.org/mathscinet-getitem?mr=#1}{#2}
}
\providecommand{\href}[2]{#2}

\end{document}